\documentclass[11pt]{article}

\usepackage{amsmath}
\usepackage{amsfonts}
\usepackage{graphicx}

\usepackage{natbib}
\bibpunct{(}{)}{;}{a}{,}{,}

\newcommand{\field}[1]{\mathbb{#1}}
\DeclareMathOperator{\PR}{\field{P}}             
\DeclareMathOperator{\E}{\field{E}}              
\def\R{\field{R}}                                

\newtheorem{theorem}{Theorem}

\newcommand{\Pf}{\paragraph{{\bf Proof.}}}       
\newcommand{\blot}{\hfill{\vrule height .9ex width .8ex depth -.1ex }}

\newcommand{\EndPf}{\ifmmode\blot\else{\unskip\nobreak\hfil
\penalty50\hskip1em\null\nobreak\hfil\blot
\parfillskip=0pt\finalhyphendemerits=0\endgraf}\fi\medskip}

\def\posRe{\R_+}
\def\Re{\R}

\def\clL{{\cal L}}

\def\eq#1/{(\ref{e:#1})}
\def\elabel#1{\label{e:#1}}
\def\flabel#1{\label{f:#1}}
\def\tlabel#1{\label{t:#1}}

\def\notes#1{}   

\def\half{\frac{1}{2}}

\def\Theorem#1{Theorem~\ref{t:#1}}
\def\Figure#1{Figure~\ref{f:#1}}

\def\eqdef{\buildrel \Delta \over =}

\def\barc{\bar c}
\def\barx{\bar x}

\def\Expect{\E}
\def\varble{\,\cdot\,}

\def\tablabel#1{\label{tab:#1}}
\def\Table#1{Table~\ref{tab:#1}}

\DeclareMathOperator{\Cov}{Cov}
\DeclareMathOperator{\Var}{Var}

\def\Ind{\field{I}}                            

\newcommand{\be}[1]{\begin{equation}\label{#1}}
\newcommand{\ee}{\end{equation}}
\newcommand{\bea} {\begin{eqnarray}}
\newcommand{\eea} {\end{eqnarray}}
\newcommand{\beas} {\begin{eqnarray*}}
\newcommand{\eeas} {\end{eqnarray*}}
\newcommand{\DefAs}{\mbox{$\,\stackrel{\bigtriangleup}{=}\,$}}
\newcommand{\refeq}[1]{(\ref{#1})}
\newcommand{\inD}{\mbox{ $\Rightarrow$ }}
\newcommand{\real}{\mbox{I\hspace{-.8mm}R}}

\newcounter{remark}
\newcounter{problem}
\newcommand{\EndEx}{\ifmmode\blot\else{\unskip\nobreak\hfil
\penalty50\hskip1em\null\nobreak\hfil\blot
\parfillskip=0pt\finalhyphendemerits=0\endgraf}\fi\medskip}

\newcounter{example}
\newenvironment{example}{\vspace{6pt} \refstepcounter{example}
{\noindent {\em Example \theexample: }}
}{\vspace{-4pt}\EndEx\par\vspace{10pt}}

\newcommand{\ExampleCtd}[1] {\vspace{6pt} {\noindent {\em Continuation
of Example \ref{#1}: }}}

\newcommand{\assump}{{\bf (A)}}
\newcommand{\bz}{{\bar z}}
\newcommand{\bZ}{{\bar Z}}
\newcommand{\by}{{\bar y}}
\newcommand{\bY}{{\bar Y}}
\newcommand{\bX}{{\bar X}}
\newcommand{\bC}{{\bar C}}

\begin{document}
\title{Variance Reduction in\\ Simulation of Multiclass Processing Networks}
\author{%
Shane G.\ Henderson%
\thanks{Work supported in part by NSF grant DMI-0224884} \\
School of Operations Research and Industrial Engineering, Cornell University \\
Ithaca, NY 14853
\and
Sean P. Meyn%
\thanks{Work supported in part by NSF grants DMI-0224884 and ECS 940372 } \\
Department of Electrical and Computer Engineering\\
University of Illinois-Urbana/Champaign \\
Urbana, IL 61801, U.S.A.
}
\maketitle

\section*{Abstract}
We use simulation to estimate the steady-state performance of a stable multiclass queueing network. Standard estimators have been seen to perform poorly when the network is heavily loaded. We introduce two new simulation estimators. The first provides substantial variance reductions in moderately-loaded networks at very little additional computational cost. The second estimator provides substantial variance reductions in heavy traffic, again for a small additional computational cost. Both methods employ the variance reduction method of control variates, and differ in terms of how the control variates are constructed.

\section{Introduction} \label{sec: intro}
{\em Owing to a mistake in the editorial process, this paper was
  accepted for publication but never actually appeared. At the request of a friend I am
  posting it on arXiv.}

\vspace{1em}
A {\em multiclass} queueing network is a network of service stations through which multiple classes of customers move. Each customer class can have different service-time characteristics at a single service station. Multiclass queueing networks are of great interest in a large variety of applications \citep{bergal87,lav83,buzsha93,ger93} because of their tremendous modeling flexibility. Perhaps the most common reason for modeling a system using a multiclass queueing network is to try to determine a suitable operating policy for the network. An {\em operating policy} is a policy that determines which customers should be worked on at which times. For example, if there are multiple customer classes at a single service station, then which class should the station work on?

In order to make comparisons between operating policies, one must define a suitable performance measure, such as expected steady-state work in process, or expected steady-state throughput, etc. For broad classes of networks one can compute certain performance measures analytically  \citep{jac63, baschamunpal75, kel79, harwil87b, harwil90a}, or one can turn to numerical computation \citep{sch84a,neu94a,daihar92,shechedaidai02}.
In general, however, these approaches are  either infeasible, or intractable due to the high complexity of network models.  Some results are available in the form of bounds on performance measures through the construction of linear programs \citep{kumkum94, berpastsi94a, kummey96a, sch01, morkum99}.
Unfortunately, these bounds are often quite loose, and so it can be difficult to compare operating policies based on such bounds alone. It is natural then, to turn to simulation.

Given that one is going to simulate many different operating policies, it is important that any simulation return relatively accurate answers as quickly as possible. This suggests the need for variance-reduction techniques that can increase the accuracy in simulation results for a given computational budget. Another reason for desiring efficient simulation techniques is that the network under consideration is often moderately to heavily loaded, in the sense that some of the resources of the network are close to full utilization. It has been noted that, in such settings, simulation can take a tremendously long time to return precise estimates of performance \citep{whi89, asm92}.

So we are strongly motivated to seek special variance-reduction techniques for multiclass networks. In this paper we develop two such variance-reduction techniques. Both are based on the approximating martingale process method \citep{hengly02, HenThesis}, which is a specialization of the method of control variates; see, for example, \citet{lawkel00} for an introduction to control variates. This paper is an outgrowth of \citet{henmey97}. The methods introduced there and here have since seen further development in \citet{henmeytad03} and \citet{bormey03a}. The theoretical results on the order of the variance constants given here have been considerably extended in \citet{mey05,mey05b}. These papers  also describe further insights on network behavior that have been uncovered since the work presented in this paper was completed.
\notes{extra reference, and a bit more humble!}
Furthermore, \citet{kimhen04} and \citet{kimhen05} have since introduced a new family of variance reduction techniques that may lead to even greater variance reductions than those seen in this paper.

Although our presentation concentrates on the estimation of the mean steady-state number of customers (of all classes) in the system, our methods may be tailored to the steady-state estimation of any linear function of the individual customer-class populations. In particular, we can also estimate, for instance, the mean steady-state number of customers of a particular class present in the system.

In Section \ref{sec: multiclass} we describe our model of a multiclass-queueing system.  We also review {\em Poisson's equation} and explain its importance in our context. In particular, we wish to approximate the solution to Poisson's equation to construct efficient simulators.

In Section \ref{sec: quadratic}, we explore quadratic forms as approximations to the solution to Poisson's equation. Computational results are given for the resulting simulation estimator, which we term the {\em quadratic} estimator. The quadratic estimator significantly outperforms a more standard simulation estimator in lightly to moderately loaded networks. In heavily loaded networks, the difference between the performance of the two estimators closes, although the quadratic estimator still provides variance reductions that would significantly reduce the computational effort involved in exploring a class of operating policies.

In Section \ref{sec: fluid}, we explore alternative approximations to the solution to Poisson's equation based on the concept of a fluid limit. The resulting simulation estimator, the {\em fluid} estimator, yields significant variance reductions in heavily loaded networks, and modest variance reductions in less heavily loaded networks. For a given simulation run-length, it is slightly more expensive to compute than the standard estimator, so that the issue of variance reduction versus computational effort needs to be considered \citep{glywhi92}. We discuss the choice of simulation estimator for a given network in Section \ref{sec: conclusions}.

\section{Multiclass Queueing Networks} \label{sec: multiclass}

Consider a system consisting of $d$ stations (or machines) and $\ell$ classes of
customers (or jobs). Class $i$ customers require service at station $s(i)$. Upon
completion of service at station $s(i)$, a class $i$ customer becomes a
class $j$ customer with probability $R_{ij}$, and exits the system
with probability
$$R_{i0} \DefAs 1 - \sum_{j=1}^\ell R_{ij}.$$

The service times for class $i$ customers are assumed to form an i.i.d.\ (independent and identically distributed)
sequence of exponentially distributed r.v.'s (random variables) with mean
$\mu_i^{-1}$.  Class $i$ customers arrive exogenously to station $s(i)$ according
to a Poisson process with rate $\lambda_i$ (which may be zero).  The Poisson-arrival processes and the service-time processes are mutually independent. We
let $\mu$ denote the $\ell$-dimensional vector of service rates, and
$\lambda$ the $d$-dimensional vector of arrival rates. Unless otherwise stated, all vectors
are assumed to be column vectors.

We require that the {\em routing matrix} $R = (R_{ij}: 1 \le i, j
\le \ell)$ be transient, so that the following inverse exists:
\[
(I-R)^{-1} = \sum_{k=0}^\infty R^k .
\]
This  ensures that all customers that enter the
system will eventually leave, and we can also be assured
that there is a unique solution $\gamma\ge0$ to the traffic
equations,
\begin{equation}
\label{eq: matrix form}
 0 = \lambda -\gamma + R' \gamma,
\end{equation}
where $R'$ denotes the transpose of the matrix $R$. We assume throughout that $\gamma_i>0$ for all $i$.

\begin{figure}[hbt]
\includegraphics[scale=0.7]{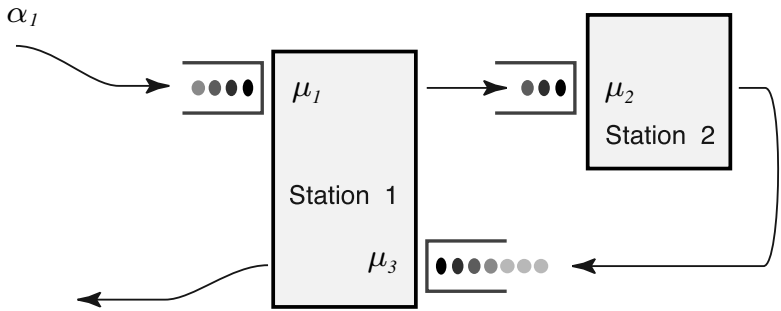}
\caption{A multiclass network with $d=2$ stations and $\ell=3$ customer classes.}
\flabel{net}
\end{figure}

\Figure{net} illustrates an example with $d=2$ stations, $\ell = 3$
customer classes, and a single exogenous arrival process (two of the
arrival rates are zero). Customer classes 1 and 3 are served at Station 1, so that $s(1) = s(3) = 1$. Similarly $s(2) = 2$. The routing values $R_{ij}$ are zero except for $R_{12} = R_{23} = 1$.

Let $X_i(t)$ denote the number of class $i$ customers present in the
system at time $t$, and let $X(t) = (X_i(t): 1 \le i \le \ell)$ be the
vector of customer populations at time $t$. Let $V_i(t)$ denote the
fraction of station $s(i)$'s effort allocated to serving customers of
class $i$ at time $t$, and let $V(t)$ denote the corresponding vector
quantity. We must have $V_i(t) \ge 0$ for all $i$ and $t$, and
$\sum_{i: s(i) = s} V_i(t) \le 1$ for all $s$ and $t$, i.e., a station $s$ can never allocate negative effort, or more
than unit effort, to the classes $\{i: s(i) = s\}$ that are served at the station.
For the network in \Figure{net} for example, suppose that at time $t$ Machine 1 is serving
Class-3 jobs, and Machine 2 is empty. Then $V(t) = (0, 0, 1)'$.

We require the operating policy adopted by stations to be {\em stationary}, {\em non-idling}, and $0-1$. We next define and explain each of these terms.

By stationary, we mean that $V(t)$ is a deterministic function of $X(t)$,
so that the workload allocations depend only on the current customer-class levels. This allows the modeling of preemptive priority policies, for example, but precludes the modeling of
policies such as FIFO which rely on additional information such as the
order in which customers arrive to a station.

By non-idling, we mean that if $\sum_{i: s(i) = s} X_i(t) > 0$, then
$\sum_{i: s(i) = s} V_i(t)=1$,
i.e., if there are customers present at station $s$ at time $t$, then
the station allocates all of its effort at that time.

By $0-1$, we mean that at any given station, at most one class receives
service at any given time. This assumption is applied {\em only} for
notational convenience.  The estimators we derive may also be applied to networks
controlled by randomized or processor-sharing policies.

With the above structure in place, we may conclude that $X = (X(t): t \ge 0)$ is a time-homogeneous
continuous-time Markov chain. However, we prefer to work in discrete
time because it simplifies the analysis. So we first rescale time so that $e' (\lambda +  \mu) = 1$, where
$e$ denotes a vector of ones. Then we uniformize; see, e.g., \citet[p.~282]{ros96}. Uniformization is a process that allows us to study the continuous-time process $X$ through a related discrete-time process $Y = (Y(n): n \ge 0)$. Define
$\tau_0 = 0$ and let the times $\{\tau_n: n \ge 1\}$ correspond to epochs
when either arrivals, real service completions, or virtual service
completions occur in the uniformized process. For $n \ge 0$, let $Y(n) =
X(\tau_n)$ and $W(n) = V(\tau_n)$. Then the process $Y= (Y(n): n \ge
0)$ is a discrete-time Markov chain evolving on a (countable) state
space $S$ that is a subset of $\{0, 1, 2, ...\}^\ell$.

\begin{example} \label{ex:mm1}
For the M/M/1 queue, $Y$ is a Markov chain on $\{0, 1, 2, \ldots\}$
with transition matrix $P$, where
$$P_{ij} = (\mu+\lambda)^{-1} \left\{
\begin{array}{ll}
\lambda & \mbox{if $j = i+1$,} \\
\mu & \mbox{if $j=\max(i-1, 0)$, and} \\
0 & \mbox{otherwise.}
\end{array}
\right.
$$
\end{example}

Our goal is to estimate $\alpha$, the steady-state mean number of
customers in the system, i.e., the steady-state mean of $|Y(0)|\DefAs e'
Y(0)$  (i.e.\  $|\varble|$ is the $L_1$ norm).   To ensure (among other
things) that
$\alpha$ exists and is finite, we make a certain assumption \assump\ below. The assumption \assump\ is known as a Lyapunov, or Foster-Lyapunov, condition. Intuitively, the function $V$ represents energy, and \refeq{eq:foster} indicates that there is an expected loss in energy for states $y$ that are ``large.'' This then ensures that energy never gets too large, and so the chain remains stable.

For any function $V$ on $ \real_+^\ell$ we define
\[
 PV(y) \DefAs \Expect_y
V(Y(1)), \qquad   \Delta_V(y) = PV\, (y) - V(y),  \qquad\qquad y\in S,
\]
where
$\Expect_y(\cdot) \DefAs \Expect(\cdot \, | \, Y(0) = y)$. Intuitively, $PV(y)$ $(\Delta_V(y))$ represents the expected energy (expected change in energy) one step from now, assuming that the chain is currently in state $y$.

\noindent \assump\ There exists a function $V: \real_+^\ell \to
\real_+$ satisfying
\begin{enumerate}
  \item $V$ is equivalent to a quintic in the sense that for some
$\delta < 1$,
\begin{equation}\label{eq:equivalent}
\delta(| y| ^5 + 1) \le V(y) \le \delta^{-1} (| y| ^5 + 1) \mbox{;
and}
\end{equation}
  \item for some $\eta > 0$, and all $y$,
\begin{equation} \label{eq:foster}
\Delta_V(y)=PV\, (y) - V(y)   \le -| y| ^4 + \eta.
\end{equation}
\end{enumerate}

\ExampleCtd{ex:mm1} For the M/M/1 queue, we may take $V(y) = b_0 y^5$,
with $b_0$ a sufficiently large constant, and then \assump\ is
satisfied as long as $\rho\DefAs\lambda / \mu < 1$.
\EndEx

In general, this condition will be satisfied if the stability-linear program of
\citet{kummey96a} admits a solution, generating a co-positive
$\ell\times\ell$ matrix $Q$.  In this case, the function $V$ may
be taken as $V(y) = (y' Q y)^{5/2}$.

Alternatively, if a fluid model is stable, and we define $\kappa =
\inf\{n \ge 1: Y(n) = 0\}$, then the function
$$V(y) = \Expect_y \sum_{k=0}^\kappa | Y(k)| ^4 $$
is bounded as in Condition 1 \citep{daimey95a}, and this function is
known to satisfy Condition 2 \citep[p.~338]{meytwe93}.

Under the assumption \assump, the chain $Y$ possesses a unique
stationary distribution $\pi$, and $\Expect_\pi | Y(0)| ^4 \le
\eta< \infty$ \citep[p.~330]{meytwe93}, where
$$\Expect_\pi(\varble) \DefAs \int_S \Expect(\varble | Y(0) = y)
\pi(dy).$$
Hence, in particular, the steady-state mean number of customers $\alpha = \Expect_\pi |Y(0)| <
\infty$.

The Lyapunov condition \assump\ is stronger than is strictly necessary
to
ensure that $\alpha = \Expect_\pi |Y(0)|$ is finite. A ``tighter''
requirement is that there is a function $h: \real_+^\ell \to
\real_+$ satisfying
\begin{equation}\label{eq:lyap}
\Delta_h(y) = Ph(y) - h(y) \le -|y| + \eta
\end{equation}
for all $y \in S$ and some constant $\eta > 0$. Given such a
function, Theorem 14.3.7 of \citet{meytwe93} allows us to conclude that
$\alpha = E_\pi |Y(0)| \le \eta$.

One might then ask whether the inequality in \refeq{eq:lyap} can be
made an equality, thereby yielding a tighter upper bound
$\eta$ on $\alpha$. In such a case we would have
\begin{equation}\label{eq:pois}
\Delta_{h^*}(y) = Ph^*(y) - h^*(y) = -|y| + \eta.
\end{equation}
This equation is known as \textit{Poisson's equation}.

If $h^*$ is a solution to Poisson's equation, then it is easy to see
that $h^* + c$ is also a solution for any constant $c$. In fact,
Proposition 17.4.1 of \citet{meytwe93} shows that any two $\pi$-integrable
solutions to Poisson's equation {\em must} differ by an additive
constant, in the sense that for all $y$ in a set $A$ with $\pi(A) = 1$,
\[
h_1(y) - \pi(h_1) = h_2(y) -\pi(h_2)
\]
where, for a real-valued function $g:S \to \real$, we denote $\pi(g) = \int_S
g(y) \pi(dy)$.

One may estimate $\alpha$ using $\alpha(n) \DefAs |\bar Y(n)|$, where
$\bar Y(n) \DefAs n^{-1} \sum_{i=0}^{n-1} Y(i)$, the mean number of customers in the system up to time $n$. The following result is a special case of Theorem 17.0.1 of \citet{meytwe93}, and shows that the
estimator $\alpha(n)$ is consistent, and satisfies a central limit theorem.

\begin{theorem} \tlabel{reptiles}
Suppose that \assump\ holds. Then $\alpha(n) \to \alpha$ almost surely (a.s.),
and furthermore,
$$n^{1/2}(\alpha(n) - \alpha) \inD \sigma N(0, 1),$$
as $n \to \infty$. The time-average variance constant (TAVC) $\sigma^2$ is
given by $\sigma^2 = \Expect[h(Y)(|Y| - \Expect|Y|)] - \Var |Y|$, where $Y$ is
distributed according to the stationary distribution $\pi$, and $h$
solves Poisson's equation; see \refeq{eq:pois}.
\end{theorem}

As noted in the introduction, it has been observed that simulation can take a very long time to yield accurate answers for heavily loaded networks \citep{asm92, whi89}. This problem is exhibited in our framework through the TAVC. Our simulation experiments indicate that the TAVC grows rapidly as the network becomes heavily loaded. Our next result lends further weight to these observations.

Consider a multiclass queueing system consisting of a single station
and (possibly) multiple customer classes. Suppose that the service
rates $\mu$ and arrival rates $\lambda$ are such that
$$e' {\cal M}^{-1} (I-R')^{-1} \lambda = 1,$$
where ${\cal M} =$ diag$(\mu)$. (This corresponds to a situation where the resources of the network are exactly matched by the demand.) Now consider a family of queueing systems
indexed by $\rho \in (0, 1)$, where the $\rho$th system has arrival
rate vector $\lambda(\rho) = \rho \lambda$. For the sake of clarity,
we occasionally suppress dependence on $\rho$.

For a given vector of buffer levels $y$, let $f(y) = d'y$
be a measure of the work in the system, where
$$d' = e' {\cal M}^{-1} (I - R')^{-1} = e'Q,$$
and $Q = {\cal M}^{-1} (I-R')^{-1}$. Intuitively, $f(y)$ measures that total expected amount of processing required to completely serve all of the customers presently in the system as given by $y$.

Let $\sigma^2_f(\rho)$ be the TAVC
associated with the estimator $f(\bY(n)) = d' \bY(n)$, and let
$\sigma^2_i(\rho)$ be the TAVC associated with $\bY_i(n)$ $(i = 1, \ldots, \ell)$. If $\Sigma_\rho$
denotes the time-average covariance matrix of $Y$, then
$\sigma^2_f(\rho) = d'\Sigma_\rho d$ and $\sigma^2_i(\rho) = e_i'
\Sigma_\rho e_i$, where $e_i$ denotes the $i$th basis vector. The proof of the following result may be found in the appendix.

\begin{theorem} \label{th:tavcisbig}
Consider the family of multiclass queueing systems above under
any non-idling work-allocation policy. Then the following are true for
each $\rho < 1$.
\begin{enumerate}
  \item Assumption \assump\ holds with the function $b_\rho f^5$ for some
sufficiently large constant $b_\rho$, and the TAVCs $\sigma^2_f(\rho)$
and $\sigma^2_i(\rho)$ $(i = 1, \ldots, \ell)$ are finite.
  \item The solution to Poisson's equation for the estimator $f(\bY(n))$ is given by
$$h(y; \rho) = \frac{f^2(y)}{2(1-\rho)} + c_\rho' y,$$
where the vector $c_\rho$ is of the order $(1-\rho)^{-1}$.
\end{enumerate}
Furthermore, there exist constants $A, B > 0$ (independent of $\rho$)
such that for $\rho$ sufficiently close to 1,
$$\frac A {(1-\rho)^4} \le \sigma^2_f(\rho) \le \frac B {(1-\rho)^4},$$
and finally
$$\mbox{trace} \,\Sigma_\rho = \sum_{i=1}^\ell \sigma^2_i(\rho) \ge \frac
{A'}{(1-\rho)^4}$$
for some constant $A'$ (again, independent of $\rho$).
\end{theorem}

Theorem \ref{th:tavcisbig} shows that trace$\,\Sigma_\rho$ is of the
order $(1-\rho)^{-4}$ as $\rho \to 1$, and this suggests, although it
does not necessarily prove, that the TAVC $e' \Sigma_\rho e$ of the
standard estimator $\alpha(n)$ is of the same order (see also  \cite{mey05} where this precise order is verified for the TAVC in a diffusion model.)  We have already mentioned that our simulation experiments also indicate that $\alpha(n)$ has high variance in heavily congested networks. Therefore,
there is strong motivation for identifying alternative estimators to
$\alpha(n)$ that can improve performance in heavy traffic. The key to
these estimators is Poisson's equation.

\medskip

If $h^*$ is $\pi$-integrable, then by taking expectations with respect
to $\pi$ in
\refeq{eq:pois}, we see that $\alpha = \eta$. Therefore, if the
solution to Poisson's equation is known, then so is $\alpha$. In
general then, we cannot expect to know the solution to Poisson's
equation. But what if an {\em approximation} is known?

If the approximation $h$ is $\pi$-integrable, then $\pi(\Delta_h)
= 0$, i.e., $\Delta_h(Y_n)$ has steady-state mean 0, and so one might consider using $\Delta_h$ in building a simulation control variate for estimating $\alpha$. In particular, one might consider using the controlled estimator
\begin{equation}\label{eq: controlled}
\alpha_c(n) = \alpha(n) + \frac \beta n \sum_{i=0}^{n-1} \Delta_h(Y(i)),
\end{equation}
where $\beta$ is an adjustable constant.

If $h = h^*$ and we take $\beta = 1$, then
$\alpha_c(n) = \alpha$, and we obtain a zero-variance
estimator of $\alpha$. In general, we can expect useful variance
reductions using the estimator $\alpha_c(n)$ provided that $h$ is a
suitable approximation to the solution to Poisson's equation. We will
discuss the choice of the constant $\beta$ later. For more
details on this approach to variance reduction see \citet{hengly02,HenThesis,henmeytad03,bormey03a} and \citet{mey05}.

So how should one go about determining an approximation to the solution
$h^*$ to Poisson's equation? It is known that for any `reasonable'
policy, the function $h^*$ is equivalent to a quadratic,  in the sense
of \refeq{eq:equivalent} (see \citet{kummey96a,mey97a,mey01a,mey05} and
\Theorem{equi} below). So it is reasonable to search for a quadratic
function $h$ that approximately solves \refeq{eq:pois}.

\section{A Quadratic Approximation} \label{sec: quadratic}

We begin this section by demonstrating the general ideas of the
approach on the stable M/M/1 queue.

\ExampleCtd{ex:mm1}
Recall that the solution to Poisson's equation is
equivalent to a quadratic. So it is reasonable to approximate the
solution $h^*$ to Poisson's equation by a quadratic function.

In the linear case $h(y) = y$, we have that $\Delta_h(y) = \lambda -
\mu w$, where $w = \Ind(y > 0)$ ($\Ind(\cdot)$ is the indicator function that is 1 if its argument is true and 0 otherwise) represents the only non-idling policy:  The
server works when customers are present, and idles when customers are
not present. Taking expectations with respect to $\pi$, we see that
the expected fraction of time that the server is working is $\Expect_\pi W(0) = \lambda / \mu$, and this is the result expected by work conservation.

Taking the pure quadratic $h(y) = a y^2$, we have
\begin{equation} \label{eq: blug}
\Delta_h(y) = 2a y(\lambda - \mu) + a \lambda  + a \mu  w.
\end{equation}
We now   choose  $a = (\mu -\lambda)^{-1} / 2$ to ensure that the
coefficient of $y$ is $-1$.  We could then use the right-hand side of
\refeq{eq: blug} as a control variate as in \refeq{eq:
controlled}. However, it is
instructive (and useful) to adopt a slightly different approach where
we avoid estimation of the variable $w$, and replace it by its known
steady-state mean $\Expect_\pi W(0)=\pi(w)=\lambda/\mu$.  We may then
use the controlled estimator $\bar Y(n) + \beta (-\bar Y(n) + 2\lambda
a)$.

This estimator is in fact equal to the estimator $\bar Y(n)+\beta\Delta_{h^*}(\bar
Y(n))$, with $h^*(y) = a(y^2 + y)$, the solution to Poisson's equation.  Hence if
$\beta = 1$, this estimator has zero variance.
\EndEx

We turn now to the general case. We assume throughout this section
that the assumption \assump\ is in place, so that the network is
stable, and all steady-state expectations that we use exist. A similar
development for reentrant lines may be found in \citet{henmey97}.

Since the solution to Poisson's equation is equivalent to a quadratic,
it is reasonable to use $h(y) = y' Q y + q'y$, for some symmetric
matrix $Q$ and vector $q$.    This approach was originally proposed in \citet{kummey96a}, based on the
prior approaches to bounding network performance presented in
\citet{kumkum94,berpastsi94a} and \citet{meydow94a}.
\notes{kumkum actually credit meydow94a for inspiration -spm}

We first consider the linear part of the
function, and then turn to the quadratic terms. Consider the
function $h_j(y) = y_j$ for some $j$. Then, letting $w$ denote the
work allocation vector corresponding to the vector $y$, we have
\begin{equation}\label{eq:itchynose}
\Delta_{h_j}(y) = \lambda_j - \mu_j w_j + \sum_i \mu_i w_i R_{ij}.
\end{equation}
Denote by $\barx$ the vector $(\mu_i \Expect_\pi  W_i(0): 1 \le i \le
\ell)$. Taking
expectations with respect to $\pi$ in \refeq{eq:itchynose} and
writing the equations (one for each $j$) in vector form, we obtain
\[
0 = \lambda -\barx + R' \barx.
\]
Since the solution to \refeq{eq: matrix form} is unique, we conclude
that
$\Expect_\pi W_j(0) = \gamma_j / \mu_j$.  Thus by considering linear
functions  we have shown that the usual traffic conditions hold.

Now consider the function $h_{jk}(y) = y_j y_k$ for $j \ne k$. Then,
\begin{eqnarray}
\Delta_{h_{jk}}(y) & = & \lambda_j y_k + \lambda_k y_j - \mu_j w_j
y_k - \mu_k w_k y_j \nonumber \\
 & + & \sum_i \mu_i w_i (R_{ij} y_k + R_{ik} y_j) -
\mu_j w_j R_{jk} - \mu_k w_k R_{kj}. \label{eq: nzpie}
\end{eqnarray}

Notice that \refeq{eq: nzpie} is a nonlinear expression in $y$,
due to the presence of the $w$ terms. We would prefer to work with
linear expressions. To this end, introduce the variables $Z_{ij}(n) =
W_i(n)Y_j(n)$, and let $\bz_{ij} = \Expect_\pi W_i(0) Y_j(0)$ for $1 \le
i, j \le \ell$.  Under our assumptions on the policy it follows that
$Z(n)$ is a fixed, deterministic function of
$Y(n)$. Furthermore, let $\by_j = \Expect_\pi Y_j(0)$.
Taking expectations with respect to $\pi$ in \refeq{eq: nzpie}, we
obtain
\begin{equation} \label{eq: balloonsup}
0 = \lambda_j \by_k + \lambda_k \by_j - \mu_j \bz_{jk} -\mu_k
\bz_{kj} + \sum_i \mu_i (R_{ij} \bz_{ik} + R_{ik} \bz_{ij}) - \gamma_j
R_{jk} - \gamma_k R_{kj}.
\end{equation}

Now, the non-idling condition implies that whenever $y_j > 0$ so that work is present at station $s(j)$,
$$\sum_{i: s(i) = s(j)} w_i = 1,$$
so that station $j$ allocates all of its effort. Hence, for {\em any} value of $y_j$, including 0,
$$y_j = \sum_{i:s(i) = s(j)} w_i y_j,$$
and consequently,
\begin{equation} \label{eq: rainagain}
\by_j = \sum_{i:s(i) = s(j)} \bz_{ij}.
\end{equation}
Therefore, if we let $\bz$ be a column vector containing the
$\bz_{jk}$'s, then the expression \refeq{eq: balloonsup} may be written
as $u_{jk}' \bz = c_{jk}$, for a suitably defined column vector $u_{jk}$ and constant $c_{jk}$.

By considering the function $h_{jj}(y) = y_j^2$, we obtain
\[
0 = 2 \gamma_j + 2 \lambda_j \by_j - 2 \mu_j \bz_{jj} + 2 \sum_i
\mu_i R_{ij} \bz_{ij},
\]
for $j = 1, \ldots, \ell$, and again these equations can be written as
$u_{jj}' \bz = c_{jj}$ for suitably defined $u_{jj}$ and $c_{jj}$.
We obtain one
equation for each $1 \le j \le k \le \ell$, so that in all there are
$\ell(\ell+1) / 2$ equations of the form $u_{jk}' \bz = c_{jk}$.
If we now write the vectors $\{ u_{jk} \}$ as columns in a matrix
$U$, and the values $\{ c_{jk} \}$ in a vector $c$, these equations can be
written as $U' \bz = c$. The matrix $U$ has $\ell^2$ rows, and
$\ell(\ell+1) / 2$ columns.

Although we began with expressions involving the Markov chain $Y$, we
are now working with the Markov chain $Z = (Z(n): n \ge
0)$. Therefore, it is useful to express the function $|y|$ as $p'
z$, for some vector $p$. In particular, $p_{ij} = 1$ if $s(i) = s(j)$,
and 0 otherwise. In view of \refeq{eq: rainagain}, the estimator
$\alpha(n)$ may be written as $p'
\bZ(n)$, where $\bZ(n) = n^{-1} \sum_{i=0}^{n-1} Z(i)$.

So define the {\em quadratic} estimator as
\begin{eqnarray}
\alpha_q(n) & \DefAs & p' \bZ(n) + \beta \nu' (U' \bZ(n) - c)
\nonumber \\
 & = & (p + \beta U \nu)' \bZ(n) - \beta \nu' c, \label{eq: hannah}
\end{eqnarray}
where $\nu$ is a vector of coefficients.

This is again of the form
\[
\alpha_q(n) = |\bar Y(n)| +\beta\Delta_h (\bar Y(n) ),
\]
where $h$ is a quadratic, $h(y) = y'Qy + \zeta'y$ for some matrix $Q$ and vector
$\zeta$.  However in the case of networks, we can only hope that $h$
\textit{approximately} solves Poisson's equation, in the sense that $Ph(y) - h(y) = - d(y)
+ \alpha_d$, with $d(\varble)\approx |\varble|$.

Let us assume (for now) that $\beta = 1$.
The variance of \refeq{eq: hannah} is then
given by
\begin{equation} \label{eq: norain}
(p + U \nu)' \Lambda_n (p + U\nu),
\end{equation}
where $\Lambda_n$ is the covariance matrix of $\bZ(n)$. But
under appropriate initial conditions and assuming \assump\ holds, $\Lambda_n \sim \Lambda / n$, where
$\Lambda$ is the time-average covariance matrix for $\bZ(n)$.
Then \refeq{eq: norain} is asymptotically given by
\begin{equation} \label{eq: speakers}
n^{-1} \| p + U \nu\|_\Lambda^2 \DefAs n^{-1} (p + U \nu)' \Lambda (p +
U\nu).
\end{equation}

Standard control variate methodology
suggests that one could estimate the covariance matrix $\Lambda$ (or
$U' \Lambda U$), and then choose $\nu$ to minimize the vector norm
\refeq{eq: speakers} \citep[p.~609]{lawkel00}. However, as cautioned in
\citet{lawkel00}, there is the danger of a variance increase associated with the
additional estimation of the covariance matrix. The ``loss factor'' is discussed
in \citet{lavwel81} and \citet{nel90} for terminating simulations, and
in \citet{loh94} for steady-state
simulation. The loss factor can become an issue when many control
variables are used (as is potentially the case here).

Rather than attempt to determine an optimal selection of control
variables from those at our disposal, we instead choose to
avoid the issue altogether by preselecting $\nu$, and then using a
standard approach to select the single
parameter $\beta$. See \citet{henmey97} for further discussion related
to this point.

From \refeq{eq: speakers}, it is ``optimal'' to choose $\nu$ to
minimize $\| U \nu + p\|_\Lambda$. Since $\Lambda$ is unknown prior to
the simulation, we instead choose $\nu$ to minimize $\| U\nu + p\|_i$
for some $L_i$ norm $\|\cdot\|$. This problem can be solved using linear
programming if $i$ is chosen to be 1 or $\infty$, or using
least-squares methods if $i = 2$.

In addition, for some workload policies, for example preemptive
priority policies, it is known that $\bz_{ij} = 0$ for some $i$ and
$j$. In this case, we may modify the norm used in the minimization
slightly to ignore the cost coefficient of $\bz_{ij}$. See
\citet{henmey97} for further discussion of this point, and
\citet{kummey96a} for the related concept of {\em auxiliary
constraints}.

It remains to explain how $\beta$ is selected. Given a response $X$,
a control variable $C$, and the form of the controlled estimator $X +
\beta C$, it is well known that the value of $\beta$ that minimizes
the controlled variance is $\beta^* = -\Cov(X, C) / \Var C$. In our
case, the response $X$ is $p' \bZ(n)$ and the control $C$ is $\nu'(U'
\bZ(n) - c)$. One may use any reasonable approach to estimate $\beta^*$
from the simulation. In particular, \citet{loh94} discusses how this
may be done in a steady-state context using both regenerative and
batch-means approaches. The process $Z$ is regenerative with
regeneration times defined by the hitting times of the state $0$, but
the regenerative cycles can be expected to be very
long. So we suggest instead using a batch-means approach to estimating
$\beta^*$. The required calculations are taken from \citet{loh94}, and
are summarized in the appendix.

\Theorem{loh}, also in the appendix, gives the relevant asymptotic theory for the
quadratic estimator. Basically, under the assumption \assump, the
quadratic estimator converges in probability, and is asymptotically
$t$-distributed when suitably normalized. The convergence mode being
``in
probability'' results from the fact that the estimator for $\beta^*$ is
only weakly consistent. If however, a strongly consistent estimator
for $\beta$ were used, or if $\beta$ were chosen to be a constant
(e.g., 1), then the quadratic estimator would be strongly consistent.

In summary, to estimate $\alpha$:
\begin{enumerate}
\item Choose $\nu$ to minimize $\| U \nu + p\|$ for some suitable norm.

\item Simulate the Markov chain $Z$ up until time $n$. (This amounts
to simulating $Y$ since $Z$ is a deterministic function of $Y$.)
\item Compute $\beta(n)$, the estimate of $\beta^*$.
\item Compute the estimator $\alpha_q(n) = p' \bZ(n) + \beta(n)
\nu'(U' \bZ(n) - c)$.
\end{enumerate}

Simulation results for the quadratic estimator and three separate queueing
networks are given in \citet{henmey97}. Simulation results
for the network in \Figure{net} operating under the FBFS (first
buffer, first-served) preemptive priority policy are given in
\Table{QuadraticResults}. We minimized $\|U \nu + p\|_2$ to
determine $\nu$,  and ignored the additional
information that $\bz_{31} = 0$ under the chosen service policy.
The results are representative of all the other networks we
experimented with.

\begin{table}[h]
\caption{Simulation Results for the Two-Station Three-Buffer
Example (to two significant figures).\tablabel{QuadraticResults}}

\begin{center}
\begin{tabular}{|l|c|c|c|c|c|}
\hline
 & \multicolumn{2}{c|}{Standard} & \multicolumn{2}{c|}{Quadratic} & \\
\hline
$\rho_2$ & Mean & Var & Mean & Var & Reduction\\
\hline
0.2 & 0.48 & 2.1E-4 & 0.48 & 1.7E-6   & 120 \\
0.4 & 1.26 & 1.4E-3 & 1.26 & 2.7E-5    & 52 \\
0.6 & 2.8  & 1.1E-2 & 2.8  & 5.1E-4     & 22 \\
0.8 & 6.9  & 0.19   & 6.9  & 2.6E-2  & 7.1 \\
0.9 & 14   & 2.0    & 14   & 0.64 & 3.1 \\
0.95 & 25  & 13     & 25   & 7.1 & 1.9 \\
0.99 & 70  & 99 & 70   & 95 & 1.0 \\
\hline
\end{tabular}
\end{center}
\end{table}

We took $\mu_1 = \mu_3 = 22$, $\mu_2 =
10$, and then chose $\lambda$ to ensure that $\rho_2 \DefAs \lambda/
\mu_2$ was as specified in the table. Before conducting the
experiments, time was rescaled so that $\lambda + \sum_i \mu_i = 1$.
Results for the standard estimator are provided in the first two
columns, and those for the quadratic estimator are in the next two
columns. The simulations were run for 100,000 time steps using 20
batches, and we repeated the experiments 200 times to obtain an
estimate of the error in the estimators. For each estimator we supply
the estimated mean and variance (over the 200 runs). The column
labeled ``Reduction'' gives the ratio of the observed variances.
The results for $\rho = 0.99$ are subject to some suspicion, owing to
the fact that our batch means exhibited correlation.

We see that for lightly to moderately loaded systems,
variance reduction factors on the order of between 120 and 7 are
observed. The largest variance reduction occurs in light traffic,
while smaller
variance reductions are obtained in moderately loaded systems.
These variance reductions are certainly useful, and come at
very little additional computational cost, since the only real
additional computational cost in computing the estimator $\alpha_q(n)$
as opposed to the standard estimator $\alpha(n)$ is the solution of
the optimization problem to choose $\nu$. But this problem is solved
once only before the simulation begins, and takes a (very)
small amount of time to solve relative to the computational effort
devoted to the simulation.

As discussed earlier, we are more interested in the performance of
these simulation estimators in moderately to heavily loaded
systems. The results in \Table{QuadraticResults} suggest that our
estimator is less effective (relative to the standard
estimator) in heavy traffic. In particular, for (very) heavily loaded
systems, only modest variance reductions are seen.

It is conceivable that the disappointing performance of the controlled
estimator in heavy traffic is due to the heuristic of choosing the
multipliers $\nu$ prior to the simulation, as opposed to attempting to
make an ``optimal'' choice. To test this idea, we estimated
the maximum possible variance reduction using the estimator
$\alpha_q(n)$. For several networks, and for various
traffic loadings, we estimated the time-average covariance matrix
$\Lambda$. As discussed earlier, the maximum possible variance
reduction (asymptotically) is obtained by selecting $\nu$ to
minimize $\| U\nu + p\|_\Lambda$. The results for the network in
\Figure{net} are presented in \Table{BestPossible}.

\begin{table}[h]
\caption{The best possible performance for the quadratic estimator in
the two-station three-buffer example (2 significant figures).
\tablabel{BestPossible}}
\begin{center}
\begin{tabular}{|l|c|c|c|}
\hline
$\rho_2$ & Standard & Quadratic & Reduction \\
\hline
0.2 & 14     & 5.7E-3 & 2500 \\
0.4 & 88     & 0.13   & 680 \\
0.6 & 670    & 2.0    & 340 \\
0.8 & 1.7E+4 & 41     & 420 \\
0.9 & 7.1E+4 & 320    & 220 \\
\hline
\end{tabular}
\end{center}
\end{table}

The columns headed ``Standard'' and ``Quadratic'' give an estimate of
the asymptotic variance for the standard and best possible quadratic
estimators respectively. The column headed ``Reduction'' gives the
ratio of these two values. The
simulation run lengths required to get reasonable estimates of $\Lambda$ for
$\rho_2 > 0.9$ were infeasibly large, and so we omitted these values.

Comparing these results with those of \Table{QuadraticResults}, we
see that the potential variance reductions appear to decrease as
congestion increases. However, substantial variance reductions may yet
be possible with a carefully chosen weighting vector $\nu$. In view of
the ``loss factor,'' the question of how best to
choose $\nu$ to achieve greater variance reduction using the quadratic
estimator is an interesting open question. In the absence of a more-effective
candidate than the one we have suggested, it would seem that a
different approach to generating control variates is warranted for
heavily loaded systems.

\section{Fluid Models and Stability} \label{sec: fluid}

To motivate our second approach to generating control variates for
the simulation of multiclass queueing networks, we consider an
alternative expression \citep[p.~432]{meytwe93} for the solution $h^*$ to
Poisson's equation \refeq{eq:pois}, namely
\begin{equation} \label{eq: exact}
h^*(y) = \Expect_y \sum_{k=0}^{\kappa} (|Y(k)| - \alpha),
\end{equation}
where $\kappa = \inf\{n \ge 0: Y(n) = 0\}$.

Although it is difficult to compute \refeq{eq: exact} exactly, we can
certainly approximate it through the use of fluid models.
As before, we introduce the key ideas through the M/M/1 queue.

\ExampleCtd{ex:mm1}
Consider a uniformized (and discrete-time) process $Y =
(Y(n): n \ge 0)$ describing the queue length  in an M/M/1 queue.
We have the recursion
$$Y(n+1) = (Y(n) + I(n+1))^+,$$
for $n \ge 0$, where $I = (I(n): n \ge 1)$ is a Bernoulli, i.i.d.\
process: $\lambda = P(I(n) = 1)$ is the arrival rate, and $\mu = P(I(n)
= -1)$ is the service rate. Time has been normalized so that $\lambda
+ \mu = 1$.

\begin{figure}[h]
\includegraphics[scale=0.5]{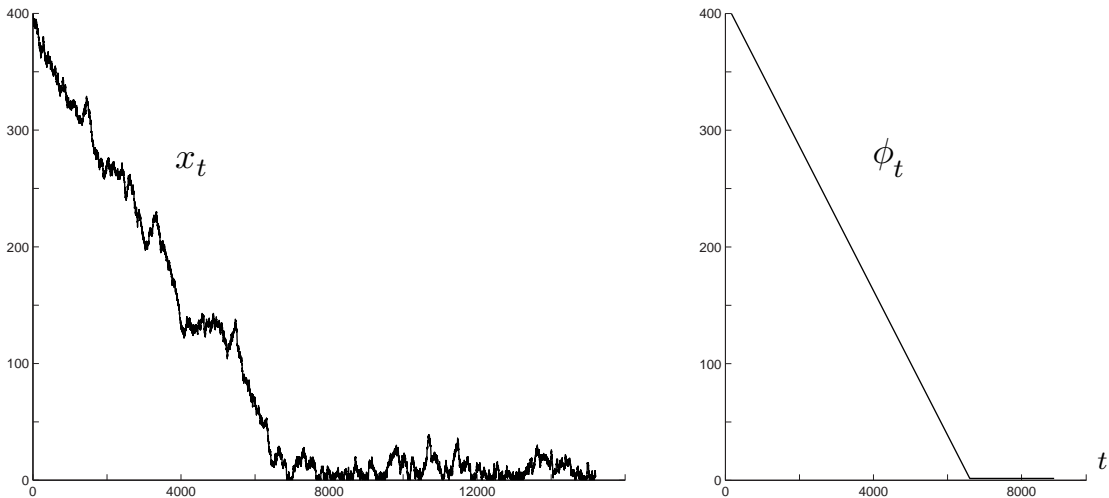}
\caption{(a) A sample path $Y$ of the M/M/1 queue with
$\rho=\lambda/\mu = 0.9$, and $Y(0)=400$. (b) A solution
to the differential equation $\dot\phi=\Ind(\phi > 0)(\lambda-\mu)$
starting from the same initial condition.\flabel{mm1} }
\end{figure}

To construct an approximation to the solution to Poisson's equation,
first note that in heavy traffic, the network will typically be
somewhat congested, and so we are primarily concerned with
``large'' states. So it may pay to consider the process starting
from a large initial condition.  In the left-hand side of \Figure{mm1}
we see one such simulation.

One approach to computing the solution to Poisson's equation is to
compute \refeq{eq: exact}.
While this is easy for the M/M/1 queue, such computation can be
formidable for more complex network models. However, consider the
right-hand side of \Figure{mm1} which shows a sample path of the
deterministic fluid, or leaky bucket model.  This satisfies the
differential equation $\dot\phi=\Ind(\phi > 0)(\lambda - \mu)$, where $\Ind(\cdot)$ is the indicator function that is 1 if its argument is true and 0 otherwise. The behavior
of the two processes looks similar when viewed on this large
spatial/temporal scale. It appears that a good approximation is
$h^*(x)\approx$
\begin{equation}
\elabel{fluidValue}
\begin{array}{rcl}
\displaystyle
h(x) &\eqdef&\displaystyle\int_0^\infty \phi(t)\, dt , \quad \phi(0)=x,
\strut
\\
\displaystyle
& =&\displaystyle \half  \frac{x^2}{\mu-\lambda}.
\end{array}
\end{equation}
This is the same approximation arrived at in Section \ref{sec:
quadratic}.
\EndEx

Of course, the M/M/1 queue is a very special case of a multiclass
queueing network, and so it is worthwhile investigating this
approximation more carefully before adopting it wholesale. We return
now to the case of a general multiclass queueing network.

The dynamics of the process $Y$ can be described by a random linear system after a slight extension of the
previous definitions. Define a sequence of i.i.d.\ random matrices $\{I(n): n \ge 1\}$ on $\{0,1\}^{(\ell + 1)^2}$, with $\PR \{\sum_j \sum_k I_{j, k}(n)  = 1\} =1$, and $\E[I_{j, k}(n)]=\mu_j R_{jk}$, where $\mu_j$ denotes the service rate for class $j$ customers. Note that exactly one element of $I(n)$ is positive for each $n$. These random variables indicate which event in the uniformized process $Y = (Y(n): n \ge 0)$ is to occur. The variable $I_{jk}(n) = 1$ if and only if a class $j$ job completes service and moves to station $k$. It is convenient to capture the exogenous arrival processes within the same framework. An exogenous arrival is indicated by $j = 0$, and a departure from the system is indicated by $k = 0$. For $j=0$, let $\mu_0 \eqdef \sum_{k=1}^\ell \lambda_k$ denote the overall arrival rate of customers to the system. Define $R_{0, 0} = 0$, and for $1 \le k \le \ell$, $R_{0, k} = \lambda_k / \mu_0$. Thus, we pool all customer arrivals into one stream with rate $\mu_0$, and an arriving customer is allocated to one of the $\ell$ classes according to the appropriate probability.

For $1 \le j \le \ell$, let $W_j(n) = 1$ if station $s(j)$ is allocating its entire effort to customers of class $j$ at time $n$, and 0 otherwise. As before we require that $W_j(n)$ is a deterministic function of $Y(n)$ for all $n \ge 0$ and all $j = 1, \ldots, \ell$. We define $W_0(n) = 1$ for all $n$, indicating that the exogenous arrival process is always active. For $1\le k \le \ell$, let $e^k$ denote the $k$th basis vector in $\Re^\ell$, and set $e^0 \eqdef 0$. The random linear system can then be defined as
\begin{equation}
 Y(n+1) = Y(n) +
 \sum_{j=0}^\ell \sum_{k=0}^\ell I_{j,k}(n+1) [-e^j + e^k ] W_j(n),
 \elabel{Net}
\end{equation}
where the state process $Y$ denotes the vector of customer classes in the system as before.

To define the fluid model associated with this network we suppose that the initial condition is large so that $m=|Y(0)|\gg 1$. We then construct a continuous time process $\phi^y (t)$ as follows:  If $tm$ is an integer, we set
\[
\phi^y (t) = \frac{1 }{ m} Y (m t),
\]
where $Y(0) = y$ and $|y| = m$. For all other $t \geq 0$, we define $\phi^y(t)$ by
linear interpolation, so that it is continuous and
piecewise linear in $t$.  Note that $|\phi^y(0)|=1$, and that
$\phi^y$ is Lipschitz continuous (see, e.g., \citet[p.~229]{apo69} for a definition
of Lipschitz continuity.) The collection of all ``fluid limits'' is defined by
\[
\clL \eqdef   \bigcap_{m=1}^\infty \overline{ \{\phi^y : |y| > m \} }
\]
where  the overbar denotes weak closure. The set $\clL$ depends upon
the particular policy chosen, and for many policies such as preemptive
priority policies, it is a family of purely deterministic functions.

Any process $\phi \in \clL$ evolves on the state space $\posRe^\ell$
and, for a wide class of scheduling policies, satisfies a differential
equation of the form
\begin{equation}
\frac{d}{dt} \phi(t) = \sum_{j=0}^\ell \sum_{k=0}^\ell  \mu_j R_{jk}
[-e^j+e^k] u_j(t)
 \elabel{fluidNet}
\end{equation}
where the function $u(\cdot)$ is analogous to the discrete control,
and satisfies similar constraints (see the M/M/1 queue model
described earlier, or \citet{dai95a} and \citet{daiwei94a} for more general
examples). In many cases the differential equation \eq fluidNet/ admits a unique
solution, from any initial condition, even though typically in practice the
control $u$ is a discontinuous function of the state $\phi$ (consider again any
priority policy).

It is now known that stability of \eq Net/ is closely connected with
the stability of the fluid
model~\citep{dai95a,kummey96a,daimey95a}. The fluid model $\clL$ is
called $L_p$-{\em stable} if
\[
\lim_{t \to \infty} \ \sup_{\phi \in \clL} \ \E  [ | \phi (t) |^p]  = 0.
\]
Let $T_0$ denote the first hitting time $\inf\{t\ge0 : \phi(t) =0\}$.
It is shown in \citet{mey97a} that
$\sup_{\phi\in\clL}\Expect[T_0]<\infty$ when the model is
$L_2$-stable. Hence, when $\clL$ is non-random, $L_2$-stability is
equivalent to \textit{stability} in the sense of \citet{dai95a}:  There
is some time $T$ such that $\phi(t)=0$ for $t\ge T$, $\phi\in\clL$. For example, in the M/M/1
queue with $\lambda < \mu$, the queue eventually hits 0 as seen in \Figure{mm1}.

The following result is a minor generalization of results from
\citet{kummey96a, daimey95a}. Its proof is omitted.

\begin{theorem}
\tlabel{equi}
The following two stability criteria are equivalent for the network
under any non-idling policy, and any $p\ge 2$.
\begin{description}
\item{(i)}
There is a function
$V$, and a constant $b<\infty$ satisfying
\[
P V(y) - V(y)  \le  -|y|^{p-1} + b
\]
where for some $\delta>0$,
\begin{equation}
 \delta(1+ |y|^p)  \le  V(y) \le  \delta ^{-1}(1+ |y|^p),
\qquad y\in S.
\elabel{QuadEquiv}
\end{equation}
\item{(ii)}
The fluid model $\clL$ is $L_p$-stable.
\end{description}
\EndPf
\end{theorem}

Thus, $L_p$ stability can be verified through the Lyapunov condition \assump. Using this result it is possible to show that the solution to Poisson's equation is asymptotically equal to a value function for the associated fluid model, provided that the fluid model is $L_2$-stable. It can be shown that many policies for the fluid model are piecewise
constant on a finite set of cones in $\real_+^\ell$.  This is certainly the case for buffer
priority policies, and also holds for $L_1$ optimal policies (for a
discussion see \citet{wei95a}). It then follows that for such policies the fluid value function $V$
is piecewise quadratic. The proof of the following result appears in the appendix.

\begin{theorem} \tlabel{h=V}
Suppose that for a given non-idling policy $w$, the fluid model $\clL$
is $L_2$-stable and non-random. Suppose moreover that limits are
unique, in the sense that $\phi_1(0)\neq \phi_2(0)$ for any two
distinct $\phi_i\in\clL$.

Then a solution $h^*$ to Poisson's equation exists, and
\[
\limsup_{ |y| \to \infty} \Bigl| \frac{h^*(y)}{ V(y)} - 1\Bigr| =0,
\]
where
\[
V(y) = |y|^2 \int_0^\infty |\phi(t)| \, dt, \qquad \phi(0)
                                                     = \frac{y}{|y|}.
\]
\end{theorem}

Hence, under the conditions of \Theorem{h=V}, a solution $h$ to Poisson's equation is intimately related to the fluid value function $V$. This result then strongly motivates the use of a fluid value function as an approximation for the solution to Poisson's equation.

Another way to motivate the fluid approximation is to note that (from \eq Net/)
\begin{eqnarray}
 Y(n+1) &=& Y(n) + \sum_{j, k} \mu_j R_{jk}(-e^j + e^k) W_j(n)
\nonumber \\
 & + &  \sum_{j, k} (I_{jk}(n+1) - \mu_j R_{jk})(-e^j + e^k) W_j(n)
\nonumber \\
 & = & Y(n) + B W(n) + D(n+1) \nonumber \\
 & = & Y(0) + \sum_{i=0}^n B W(i) + M(n+1), \label{eq:Net2}
\end{eqnarray}
where $B$ is an $(\ell+1) \times
(\ell + 1)$ matrix. The process $M(\cdot)$ is a vector-valued
martingale with respect to the natural filtration, and $D(i)$ is the
martingale difference $M(i) - M(i-1)$. (See, e.g., \citealt{ros96} for an introduction
to martingales.) It is straightforward to check
that $E D(i)' D(i)$ is bounded in $i$ (by $b$ say), so that $E M(n)'
M(n)
\le b n$ for all $n \ge 0$. Hence, the network is essentially a
deterministic fluid model with a `disturbance' $M$.  When the initial
condition $Y(0)$ is large, then the state dominates this disturbance,
and hence the network behavior appears deterministic.

We therefore have strong motivation for approximating the solution to
Poisson's equation $h^*$ by $V$, where
$$V(y) \DefAs \int_0^\infty \phi(t)\, dt,$$
and $\phi$ solves the differential equation \eq fluidNet/ with
$\phi(0) = y$.
The {\em fluid estimator} of $\alpha$ is then
\begin{equation} \label{eq: fluid est}
\alpha_f(n) = |\bY(n)| + \frac \beta n \sum_{k=0}^{n-1}\Delta_V(Y(k)).
\end{equation}
The parameter $\beta$ is again a constant that may be chosen to
attempt to minimize the variance of the fluid estimator. We use the
methodology outlined in the appendix to estimate the optimal
$\beta$. Consequently, the asymptotic results for the quadratic
estimator also apply here, namely that under the assumption \assump,
the fluid estimator is weakly consistent, and is asymptotically
$t$-distributed when suitably normalized.

Clearly, to implement the fluid estimator we need to be able to compute
$\Delta_V$. For any function $V$, we have that
\begin{eqnarray*}
PV(y) & = & \Expect_y V(Y(1)) \\
 & = & \sum_{j=0}^\ell \sum_{k=0}^\ell \mu_j R_{jk} V(y - e^j +
e^k) W_j(y),
\end{eqnarray*}
so that it suffices to be able to compute $V(y)$ for $y \in S$.

Solving the differential equation \eq fluidNet/ to find $\phi$ is not difficult when the fluid control $u$ is
piecewise constant on a finite set of cones in $\real_+^\ell$
since in this case $\phi$ is piecewise linear. Integrating $\phi$ to find the fluid value
function $V$ is then straightforward. In other words, for a specific model, some preliminary work has to be done to give code that can compute $V$, but this is usually not a difficult step.
Algorithms for computation or approximation of $V$ are described in \citet{enghummey96a} and \citet{bormey03a}.

It should be apparent from the above discussion that
computing the control $\Delta_V$ for the estimator \refeq{eq: fluid est}
may be moderately time-consuming (computationally speaking) relative to
the time taken to simply simulate the process $Y$. Be that as it may, it
is certainly the case that the time taken to compute the control is
relatively insensitive to the congestion in the system.

In \Table{FluidResults} we present simulation results for the fluid
estimator on the network
of \Figure{net}. (Similar results were obtained for all other
multiclass queueing networks that we tried.) The entries in
\Table{FluidResults} have the same interpretation as those in
\Table{QuadraticResults}. In particular, the column headed
``Reduction'' represents the variance reduction factor over the
standard estimator.

\begin{table}[h]
\caption{Simulation results for the two-station three-buffer
example (2 significant figures). The interpretation of the values given is the same as in \Table{QuadraticResults}. \tablabel{FluidResults}}

\begin{center}
\begin{tabular}{|l|c|c|c|c|c|}
\hline
$\rho_2$ & Mean & Var & Reduction\\
\hline
0.2 & 0.47 & 6.1E-5 & 3.5 \\
0.4 & 1.3  & 4.0E-4 & 3.5 \\
0.6 & 2.8  & 3.5E-3 & 3.1 \\
0.8 & 6.9  & 4.3E-2 & 4.4 \\
0.9 & 14   & 0.17   & 12  \\
0.95 & 26  & 0.23   & 56  \\
0.99 & 110 & 0.98   & 100 \\
\hline
\end{tabular}
\end{center}
\end{table}

The best value of $\beta$ was found to be close to unity in each of the
simulations, particularly at high loads where it was found to be within $\pm
5$\% of unity.

Observe that for low traffic intensities, the fluid estimator yields
reasonable variance reductions over the standard estimator. However,
because it is more expensive to compute than the standard estimator,
these results are not particularly encouraging. But as the system
becomes more and more congested, the fluid estimator yields large
variance reductions over the standard estimator, meaning that the
extra computational effort per iteration is certainly worthwhile. For
very high traffic intensities, the fluid estimator {\em
significantly} outperforms both the standard estimator and the
quadratic estimator, and so we have achieved our goal of deriving an
estimator that can be effective in heavy traffic.

\section{Conclusions} \label{sec: conclusions}

We have given two simulation estimators for estimating a linear
function of the steady-state customer class population.

The quadratic estimator produces very useful
variance reductions in light to moderate traffic at very little
additional computational cost. We recommend that it be used in
simulations of such lightly loaded networks. The quadratic
estimator is less effective in simulations of heavily loaded networks,
but could potentially provide useful variance reductions in this regime if a
better choice of weighting vector $\nu$ can be employed.

The fluid estimator provides modest variance reduction in light to
moderate traffic, but appears to be very effective in heavy traffic. There
is an additional computational overhead in computing the fluid
estimator, but this overhead is (roughly) independent of the load on
the network. Hence we may conclude that in heavily loaded systems, the
fluid estimator should yield significant computational improvements,
and should therefore be used.

One might conclude from the above discussion that the quadratic and
fluid estimators could be combined using the method of multiple
control variates to yield a single ``combined'' estimator. However, we
believe that it is unlikely that a combined estimator would yield
significant improvements over the use of either the quadratic
estimator (in light traffic) or the fluid estimator (in heavy
traffic). In light traffic, we expect that the additional reductions
in variance would be negated by the increased computational
effort. And in heavy traffic we expect that any additional variance
reduction would be modest, owing to the weaker performance of the
quadratic estimator in this regime.

Finally, we note that \citet{vea04a} explores bounded perturbations of the fluid value function in approximate dynamic programming.  Related techniques are   considered in current research to refine the fluid estimator.
\notes{note!}

\bibliographystyle{abbrvnat}

\section*{Appendix}

\subsection*{Proof of Theorem \ref{th:tavcisbig}}
It is straightforward to show that
\begin{equation} \label{eq: dos}
\Delta_{f^2}(y) = -2(1-\rho)f(y) + e'Q [L + Q^{-1}W - {\cal M}WR +
\mbox{diag}(e'{\cal M}WR) ]Q'e,
\end{equation}
where $W$ is the diagonal matrix containing the work allocation vector $w$
corresponding to $y$, and $L =$ diag$(\lambda)$. Furthermore,
for $m \ge 3$,
\begin{equation} \label{eq: cinque}
\Delta_{f^m}(y) = -m(1-\rho) f^{m-1}(y) + \mbox{lower order terms}.
\end{equation}
It follows from \refeq{eq: cinque} with $m = 5$ that \assump\ holds.
This implies (Theorem 17.5.3 of \citet{meytwe93}) that the TAVC's are finite,
and furthermore, that
$$\sigma^2_f(\rho) = \lim_{n \to \infty} n \Var(d' \bY(n)), \mbox{ and}$$
$$\sigma^2_i(\rho) = \lim_{n \to \infty} n \Var(\bY_i(n)).$$

The fact that $h$ is of the form given in the theorem follows from
(\ref{eq: dos}). Observe that the function $f^2 / (2(1-\rho))$ is
``almost'' the solution to Poisson's equation. It needs to be
adjusted slightly to remove the terms in the RHS of (\ref{eq: dos})
involving the work allocation vector $w$. These terms are of the form
$f'w$, where the coefficients in $f$ are bounded in $\rho$. Therefore,
the solution to Poisson's equation is as given in the theorem.

So then the TAVC $\sigma^2_f(\rho)$ is given by
\begin{equation} \label{eq:coldfingers}
\Expect h(Y)(f(Y) - d'\by) = \frac{\Cov(f^2(Y), f(Y))}{2(1-\rho)}
- \Expect(c_\rho'Y \, (f(Y) - d'\by)),
\end{equation}
where $Y$ is distributed according to the stationary distribution
$\pi$ and $\by = \Expect Y$.

Since $\Expect \Delta_{f^m}(Y) = 0$ for $m = 0, \ldots, 4$, it follows from
\refeq{eq: dos} that $\Expect f(Y)$ is of the order $(1-\rho)^{-1}$ as
$\rho \to 1$. Then, by induction using \refeq{eq: cinque}, $\Expect
f^m(Y)$ is of the order $(1-\rho)^{-m}$ as $\rho \to 1$ for $m = 1, \ldots, 4$.

The second term on the right-hand side of \refeq{eq:coldfingers} is
therefore of the order $(1-\rho)^{-3}$ as $\rho \to 1$. As for the
first term, we have the easily proved inequality that for any
non-negative r.v.\ $X$ with $\Expect X^3 < \infty$,
$$\Cov(X^2, X) \ge \frac{\Var X}{\Expect X^2} \Expect X^3.$$
Applying this inequality to the first term on the right-hand side of
\refeq{eq:coldfingers}, and noting that $\Var f(Y)$ is of the same
order as $\Expect f(Y)^2$ as $\rho \to 1$, we obtain the required result
that $\sigma^2_f(\rho)$ is of the order $(1-\rho)^{-4}$ as $\rho \to 1$.

The last statement of the theorem follows from the fact that
\begin{eqnarray*}
\sigma^2_f(\rho) & = & \lim_{n \to \infty} n \Var(d' \bY(n)) \\
 & = & \lim_{n \to \infty} n \Var(\sum_{i=1}^\ell d_i \bY_i(n)) \\
 & \le & \lim_{n \to \infty} \ell n \sum_{i=1}^\ell \Var(d_i \bar Y_i(n)) \\
 & = & \ell \sum_{i=1}^\ell d_i^2 \sigma^2_i(\rho).
\end{eqnarray*}
\EndPf

\subsection*{Control Variates in Steady-State Simulation}
We repeat formulae from \citet{loh94} for
estimating the control variate parameter $\beta$ that is used in both
the quadratic estimator \refeq{eq: hannah} and the fluid estimator
\refeq{eq: fluid est} via
the batch means method of simulation output analysis. To
encapsulate both estimators and avoid repetition, we give the formulae
for the case where a real-valued stochastic process $X =
(X(n): n \ge 0)$ is simulated, and a real-valued control $C = (C(n): n
\ge 0)$ is recorded.

Let the $b$ batches each consist of $m$ observations, so that the
simulation run-length $n=mb$. If
$X_i$ and $C_i$ are the $i$th batch means of the process and
control respectively, then for $0 \le i \le b-1$, we have
$$X_i = \frac 1 m \sum_{j=im}^{(i+1)m - 1} X(j) \mbox{ and }
C_i = \frac 1 m \sum_{j=im}^{(i+1)m - 1} C(j).$$
Let $\bX_n$ and $\bC_n$ denote the overall (sample) means of the
process and control respectively.

Define
\begin{eqnarray*}
V_{XX}(n) & = & \frac 1 {b-1} \sum_{i=0}^{b-1} (X_i - \bX_n)^2, \\
V_{CC}(n) & = & \frac 1 {b-1} \sum_{i=0}^{b-1} (C_i - \bC_n)^2, \mbox{
and}\\
V_{XC}(n) & = & \frac 1 {b-1} \sum_{i=0}^{b-1} (X_i - \bX_n)(C_i -
\bC_n).
\end{eqnarray*}

Define $\beta = -V_{XC} / V_{CC}$, and let $\alpha_n = \bX_n + \beta
\bC_n$ be the controlled estimator.

Finally, let
\[
R^2(n) = \frac{b-1}{b-2}\left(V_{XX}(n) -
\frac{V_{XC}(n)^2}{V_{CC}(n)} \right)
\]
and
\begin{equation} \label{eq:ick}
S^2(n) = R^2(n) \left( \frac 1 b + \frac 1 {b-1} \frac
{\bC_n^2}{V_{CC}} \right).
\end{equation}

Using the above computational process, we can construct the
quadratic estimator $\alpha_q(n)$ and the fluid estimator $\alpha_f(n)$.
The following result describes the asymptotic behavior of these
estimators.

\begin{theorem}[Loh 1994] \tlabel{loh}
Under the assumption \assump, the estimators $\alpha_q(n)$ and
$\alpha_f(n)$ converge in probability to $\alpha$, and for $j = q, f$,
$$\frac{\alpha_j(n) - \alpha}{S_j(n)} \inD T_{b-2},$$
where $T_{b-2}$ has the Student's $t$-distribution with $b-2$ degrees
of freedom, and $S_j(n)$ is defined in the obvious way through
\refeq{eq:ick}.
\end{theorem}
\Pf
The second result is proved in Section 1.3.2 of \citet{loh94}. In
addition, Proposition 1.5 of \citet{loh94} shows that $nS_j^2(n)$
converges in distribution to a finite-valued random variable as $n\to
\infty$ so that the first result follows.
\EndPf

\subsection{Proof of \Theorem{h=V}}
First note that for any $\phi\in\clL$ we have $|\phi(0)|=1$, and under
the assumptions of the theorem there is a   $T>0$ such that
\[
V(y) = |y|^2 \int_0^T |\phi(t)| \, dt, \qquad \phi(0)=\frac{y}{|y|} ,\quad
\phi\in\clL.
\]
We shall fix such a $T$ throughout the proof.

 From \Theorem{equi}, a Lyapunov function exists that satisfies
\refeq{eq:lyap}, which is equivalent to a quadratic in the
sense of \eq QuadEquiv/.  It follows that $\pi (c) <\infty$, where
$c(y) = |y|$, and that a solution to Poisson's equation exists which is
bounded from above by a quadratic, and uniformly bounded from below
\citep[p.~432]{meytwe93}.

We can then take the solution to Poisson's equation and iterate as
follows: $P^nh = h - \sum_{i=0}^{n-1} P^i\barc$, where $\barc(y) = |y| -
\alpha$.  Let $m=|y|$ and take  $n = [m T] =$ the integer part of
$mT$ to give
\[
\frac{\Expect_{ y}[h(Y( m T))]}{ m^2}
= \frac{h( y)}{ m^2}
        - \frac{\Expect_{ y}\Bigl[ \sum_{i=0}^{ [m T]-1}
\frac{|Y(i)|}{m}\Bigr]}{m}
        - \frac{T\alpha}{ m}.
\]

Since $h$ is bounded above by a quadratic, there is a $K<\infty$ such
that
\[
\left| \frac{[h(Y( m T))]}{ m^2} \right|
 \le K\Bigl(1+ \frac { |Y( m T)|^2}{ m^2}\Bigr)
\]
The random variable on the right hand side is
uniformly bounded by $K(1+1/m+ T)^2$ for all initial $y$ since at most one
customer can arrive during each time slot.  It then follows from weak convergence
(see, e.g., \citet{bil68}) and the definition of
$T$ that
\[
\limsup_{|y|\to\infty} \Expect_y \left| \frac{[h(Y( m T))]}{ m^2}
\right|
=0
\]

Moreover, again by Lipschitz continuity of the fluid model we have
\[
\limsup_{|y|\to \infty}
\left| \Expect_{ y}\left[
\frac{1}{ m}
\sum_{i=0}^{ [m T]-1} \frac{|Y(i)|}{ m}\right] - V(\frac{y}{m})\right| =
0.
\]
Putting these results together we see that
\[
\limsup_{|y|\to \infty}
\left|  \frac{h(y)}{m^2} -
V(\frac{y}{m})\right| = 0,
\]
proving the result.
\EndPf

\end{document}